\documentclass[letterpaper, 10 pt, conference]{ieeeconf}  

\IEEEoverridecommandlockouts                          
\usepackage{cite}
\usepackage{amsmath,amssymb,amsfonts,bm}
\usepackage{bbm}
\usepackage{graphicx}
\usepackage{subfigure}
\usepackage{textcomp}
\usepackage{algorithm}
\usepackage{algpseudocode}

\usepackage{afterpage}
\usepackage{stfloats} 
\usepackage{cuted}

\usepackage{mathrsfs}
\usepackage{xcolor}
\allowdisplaybreaks[4]
\newtheorem{definition}{Definition}
\newtheorem{theorem}{Theorem}
\newtheorem{lemma}{Lemma}
\newtheorem{assum}{Assumption}
\newtheorem{coro}{Corollary}

\title{\LARGE \bf
Adaptive Control of Unknown Linear Switched Systems via Policy Gradient Methods
}

\author{Felix Laurent, Feiran Zhao, Jaap Eising, Florian D\"{o}rfler 
\thanks{This work was supported by ETH Zurich and the SNF through the NCCR Automation.}
\thanks{F. Laurent is with the Department of Mechanical and Process Engineering, ETH Zurich, 8092 Zurich, Switzerland. (e-mail: flaurent@ethz.ch)}
\thanks{F. Zhao, J. Eising, and F. D\"{o}rfler are with the Department of Information Technology and Electrical Engineering, ETH Zurich, 8092 Zurich, Switzerland. (e-mail: zhaofe@control.ee.ethz.ch, jeising@control.ee.ethz.ch, dorfler@ethz.ch)} %
}

\begin{document}

\maketitle
\thispagestyle{empty}
\pagestyle{empty}
\begin{abstract}
We consider the policy gradient adaptive control (PGAC) framework, which adaptively updates a control policy in real time, by performing data-based gradient descent steps on the linear quadratic regulator cost. This method has empirically shown to react to changing circumstances, such as model parameters, efficiently. To formalize this observation, we design a PGAC method which stabilizes linear switched systems, where both model parameters and switching time are unknown. We use sliding window data for the policy gradient estimate and show that under a dwell time condition and small dynamics variation, the policy can track the switching dynamics and ensure closed-loop stability. We perform simulations to validate our theoretical results.
\end{abstract}

\section{Introduction}

The linear quadratic regulator (LQR) is a standard benchmark in control theory, reinforcement learning, and their applications. Traditionally, the design pipeline for an LQR controller design was \textit{indirect}, first a model of the to-be-controlled system is determined, after which the controller is designed. With the recent interest in data-driven control, came contrasting \textit{direct} approaches: the aim of learning a controller based on measurements of a system. 

This problem has been viewed through a number of different lenses. In the wake of methods based on Willems' fundamental lemma~\cite{willems2005note}, the finite-horizon case has attracted a lot of attention in receding horizon formulations, starting from~\cite{coulson2019data,berberich2020data}. In contrast to these, the one-shot design of controllers for the infinite-horizon LQR is studied in ~\cite{de2019formulas,van2020data}. 
Simultaneously, the convergence of learning controllers was studied using tools from reinforcement learning~\cite{fazel2018global,mohammadi2022convergence,tsiamis2023statistical}. 
These methods generally deal with data that is collected in episodes, or even off-line: the controller is learned from batches of data during which the policy does not change, and learning and control are thus separated. In turn, this means that the aforementioned methods are unable to improve the policy efficiently in real time.

In comparison, adaptive control can monitor its own performance and adjust its parameters in the direction of better performance~\cite{drenick1957adaptive, anderson2005failures}. Classical adaptive control focuses on robust stability, designing the control input with online closed-loop data based on Lyapunov theory~\cite{ioannou1996robust}. Representative instances include model-reference adaptive control (MRAC), self-tuning regulators, and model-free adaptive control~\cite{ioannou1996robust,8621060}. There are also adaptive control methods that seek both stability and convergence of the policy using certainty-equivalence LQR~\cite{mania_certainty_2019,lu2023almost,simchowitz2020naive,cohen2019learning}. While these advances are achieved for linear time-invariant (LTI) systems, adaptation of the control policy plays an important role where systems are not LTI~\cite{dai2025online, teutsch2023online, wang2024data}. In particular, when the unknown dynamics are switching, the policy must be updated in real time with online data to achieve stability. Recently, the behavioral system theory and data-informativity approaches to data-driven control have been extended to the situation of linear switched systems~\cite{rotulo2022online, liu2023online,bianchi2025data}. However, these methods rely on solving an optimization problem at each time step, which could prove to be computationally intensive, and sensitive to noise.

Another approach to achieve adaptation is gradient-based, updating the policy iteratively with gradient descent in real time. Recently, we have proposed a policy gradient adaptive control (PGAC) framework \cite{zhao2025policy}, which adaptively updates the policy in feedback by descending the gradient of the LQR cost efficiently computed from data. It includes both direct and indirect gradient computation approaches: the indirect approach estimates a model, based on which the policy gradient is computed; the direct approach uses a sample covariance parameterization of the LQR and estimates the gradient directly from data \cite{zhao2024data}. For linear systems, the proposed policy gradient methods exhibit closed-loop stability and convergence of the policy to the optimal LQR gain under mild assumptions. Moreover, though theoretical guarantees have been shown only for LTI systems, the PGAC framework has been successfully applied in nonlinear or time-varying situations, such as power converter systems~\cite{11187193}, aerospace control~\cite{10981973}, and autonomous bicycle control~\cite{persson2025adaptive}.

In this paper, we focus on designing a PGAC method to stabilize linear switched systems. Unlike the literature \cite{sun2005switched}, we assume that both the model parameters and the switching times are unknown. We adopt the infinite-horizon LQR cost as the objective function, as the optimal LQR gain can stabilize the system and has a certain degree of stability margin. Then, we use indirect PGAC method, where the gradient of the cost is computed based on an identified dynamical model. To make this gradient reflect the varying dynamics well, we use sliding window data for model identification. We show that under a dwell time condition and a bounded variation of the dynamics, the method can track the switching dynamics and ensure closed-loop stability. 
Simulations are conducted to validate our theoretical results.

The remainder of this paper is organized as follows. Section II recapitulates the indirect PGAC method and formulates the linear switched system control problem. Section III proves the main results on system stability. Section IV performs numerical case studies. Conclusions are made in Section V. All proofs are deferred to the Appendix.

\section{Preliminaries}
In this section, we recapitulate the \textit{policy gradient adaptive control} (PGAC) algorithm \cite{zhao2025policy}.  
\subsection{Policy gradient methods for the linear quadratic regulator}
Consider the linear time-invariant (LTI) system  
\begin{equation}\label{equ:sys}
x_{t+1} = A x_t + B u_t + w_t,
\end{equation}
where $t\in \mathbb{N}$, $x_t\in\mathbb{R}^{n}$ is the state, $u_t\in\mathbb{R}^{m}$ is the control input, and $w_t\in \mathbb{N}$ is the noise. We assume that the system~\eqref{equ:sys} is controllable. It is well known that the linear quadratic regulator (LQR) cost can be regarded (see e.g. \cite{fazel2018global}) as a function of a state feedback gain $K\in \mathbb{R}^{m\times n}$: 
\begin{equation}
    C(K) := \mathrm{Tr} \left( (Q + K^\top R K) \Sigma_K \right),
\end{equation}
where the weighting matrices $(Q, R)$ are positive definite, and \( \Sigma_K \) is the closed-loop state covariance matrix obtained as the positive definite solution to the Lyapunov equation
\begin{equation}
\label{Sigma}
\Sigma_K = I_n + (A + BK) \Sigma_K (A + BK)^\top.
\end{equation} 
Instead of minimizing this cost, the policy gradient method updates the feedback gain, or \textit{policy}, using gradient descent:
 \begin{equation}\label{equ:grad}
 	 	K_{t+1} = K_t - \eta \nabla C(K_t),
 	 \end{equation}
     where $\nabla C(K_t)$ is the gradient of $C$ with respect to $K_t$. Define $\mathcal{S}:=\{K\in \mathbb{R}^{m\times n}|\rho(A+BK)<1\}$ as the set of stabilizing gains. Then, the gradient $\nabla C(K)$ admits a closed-form expression.
 	 \begin{lemma}[{\cite[Lemma 1]{fazel2018global}}]\label{lem:mb_pg}
 	 	For any $K\in\mathcal{S}$, the gradient of $C(K)$ is given by
 	 	\begin{equation}\label{equ:gradK}
 	 		\nabla C(K) =  2\left(\left(R+{B}^{\top} P {B}\right) K+{B}^{\top} P {A}\right) {\Sigma_K},
 	 	\end{equation} 
 	 	where ${\Sigma_K}$ satisfies \eqref{Sigma}, and $P$ is the positive definite solution to the Lyapunov equation 
 	 	\begin{equation}\label{equ:Lyap}
 	 		P = Q + K^{\top}RK + ({A}+{B}K)^{\top}P ({A}+{B}K).
 	 	\end{equation}
 	 \end{lemma}

In this paper, we consider the case where $(A,B)$ are unknown and we are unable to compute the gradient directly. In this case, PGAC uses online identification of the model for gradient computation, as recapitulated below.

    \subsection{Policy gradient adaptive control}
Since the model $(A,B)$ is unknown, we identify it from measured data \cite{dorfler2021certainty}. Consider $t$-long time series of states, inputs, noises and successor states
\begin{equation}
\label{time_series}
    \begin{aligned}
        X_{0,t} &:= [x_0,x_1, \ldots, x_{t-1}] &\in \mathbb{R}^{n\times t}, \\
        U_{0,t} &:= [u_0,u_1, \ldots, u_{t-1}] &\in \mathbb{R}^{m\times t}, \\
        W_{0,t} &:= [w_0,w_1, \ldots, w_{t-1}] &\in \mathbb{R}^{n\times t}, \\
        X_{1,t} &:= [x_1,x_2, \ldots, x_t] &\in \mathbb{R}^{n\times t},
    \end{aligned}
\end{equation} 
which satisfy the linear dynamics
\begin{equation}
\label{system}
    X_1 = AX_{0,t} + BU_{0,t} + W_{0,t}.
\end{equation}
Suppose that the data matrix $D_{0,t} = [U_{0,t}^{\top}, X_{0,t}^{\top}]$ has full row rank. Then, an estimated model can be obtained from the unique solution of the ordinary least-square problem 
\begin{equation}
\label{OLS}
    [\hat{B}_t, \hat{A}_t] = \underset{[B, A]}{\text{arg min}} \| X_{1,t} - [B,A]D_{{0,t}}
    \|_F = X_{1,t} D_{0,t}^\dagger,
\end{equation}  
where $(\cdot)^\dagger$ denotes the Moore-Penrose pseudoinverse. The gradient of the cost function $\hat{C}(K)$ can then be computed by using the model estimate $(\hat{A}, \hat{B})$ as ground-truth.  

By iteratively identifying the dynamics, applying a gradient descent step to the feedback gain, and using the resulting policy to gather closed-loop data online, the method enables adaptive LQR learning with optimal non-asymptotic performance guarantees \cite{wang2021exact, lu2023almost}. The approach is summarized in Algorithm \ref{algo2}. It requires an initial stabilizing gain and a batch of offline PE data. In the online stage, the control input includes a probing noise $e_t$ to ensure that the data is informative. We refer to \cite{zhao2024data,zhao2025policy} for more details.

\begin{algorithm}
\caption{Policy gradient adaptive control for LIT systems}
\label{algo2}
\begin{algorithmic}[1]
\State \textbf{Initialize:} Offline data \((X_{0,t_0}, U_{0,t_0}, X_{1,t_0})\), an initial stabilizing policy \(K_{t_0}\), and a stepsize \(\eta\).
\For{ \(t = t_0, t_0 + 1, \ldots\) }
    \State Apply \(u_t = K_t x_t + e_t\) and observe \(x_{t+1}\).
    \State Estimate model parameters via ordinary least-squares
    \[
    [\hat{B}_{t+1}, \hat{A}_{t+1}] = \mathop{\arg\min}_{B, A} \| X_{0,t+1} - [B, A] D_{0,t+1} \|_F.
    \]
    \State Perform one-step policy gradient descent
    \[
    K_{t+1} = K_t - \eta \nabla \hat{C}_{t+1}(K_t),
    \]
    where \(\nabla \hat{C}_{t+1}(K_t)\) is the policy gradient in Lemma \ref{lem:mb_pg} with the estimated model \((\hat{A}_{t+1}, \hat{B}_{t+1})\).
\EndFor
\end{algorithmic}
\end{algorithm}

\section{Policy gradient adaptive control for linear switched systems}

The PGAC method in Algorithm \ref{algo2} has been certified theoretically only for LTI systems \cite{zhao2024data,zhao2025policy}.
In this section, we extend the PGAC method for stabilization of switched systems.

\subsection{Problem formulation and the PGAC method for linear switched systems}
Consider the following switched system
\begin{equation}\label{equ:switch}
x_{t+1} = A_i x_t + B_i u_t, \quad T_{i-1} \leq t < T_{i}
\end{equation} 
for $i \in \mathbb{N}$, where $T_i$ is the switching time between the $i^{th}$ and the $(i+1)^{th}$ modes. We assume that both dynamics $(A_i, B_i)$ and $T_i$, $\forall i\in \mathbb{N}$ are unknown. Our goal is to develop a PGAC method to stabilize \eqref{equ:switch}. Clearly, it is not sufficient to stabilize the system at each time $t$, and therefore, we use the notion of sequential strong stability \cite{cohen2018online}. 



\begin{definition}[\textbf{Sequential strong stability}] A sequence of policies \( K_t \) for the system \eqref{equ:switch} is \( (\kappa, \alpha) \)-strongly stable (for \( \kappa \geq 1 \) and \( 0 < \alpha \leq 1 \)) if there exist matrices \( H_t \succ 0 \) and \( L_t \) such that
\[A_t + B_t K_t = H_t L_t H_t^{-1} \quad \text{for all } t,\]
and
\begin{itemize}
  \item[(i)] \( \|L_t\| \leq 1 - \alpha \) and \( \|K_t\| \leq \kappa \);
  \item[(ii)] \( \|H_t\| \leq \kappa \) and \( \|H_t^{-1}\| \leq 1 \);
  \item[(iii)] \( \|H_{t+1}^{-1} H_t\| \leq 1 + \alpha/2 \).
\end{itemize} 
\label{def:sss}
\end{definition}  

In particular, this implies that every policy $K_t$ must be individually stable ((i) and (ii)), and that the changes between two consecutive closed-loop systems must be small (iii). Sequential strong stability gives a quantitative definition of stability and implies asymptotic stability for the state. Our aim is now to adapt PGAC to yield sequential strong stability for the system~\eqref{equ:switch}. 

To avoid accumulating errors, we use \textit{sliding window data} for identification, that is, we use only the most recent data. Define the data series with fixed window size $L$ as 
\begin{equation}\label{equ:swd}
    \begin{aligned}
        X_t &:= [x_{t-L}, x_{t-L+1}, ...,x_{t-1}] \in \mathbb{R}^{n\times L}, \\
        U_t &:= [u_{t-L}, u_{t-L+1}, ...,u_{t-1}] \in \mathbb{R}^{m\times L}, \\
        X_{t+1} &:= [x_{t-L+1}, x_{t-L+2}, ...,x_{t}] \in \mathbb{R}^{n\times L}.
    \end{aligned}
\end{equation}
To allow that the full row rank condition of $D_t = [U_t^{\top}, X_t^{\top}]$ is satisfied, we require \(L \geq n+m\). The resulting PGAC method is similar to Algorithm \ref{algo2}; see Algorithm \ref{algo2_switch}.

\begin{algorithm}
\caption{Policy gradient adaptive control for linear switched systems with sliding window data}
\label{algo2_switch}
\begin{algorithmic}[1]
\State \textbf{Initialize:} Offline data \((X_{t_0}, U_{t_0}, X_{t_0+1})\), an initial stabilizing policy \(K_{t_0}\), and a stepsize \(\eta\).
\For{ \(t = t_0, t_0 + 1, \ldots\) }
    \State Apply \(u_t = K_t x_t + e_t\), then observe \(x_{t+1}\).
    \State Estimate a model using $(X_{t+1}, U_{t+1}, X_{t+2})$ and ordinary least-squares
    \[
    [\hat{B}_{t+1}, \hat{A}_{t+1}] = \mathop{\arg\min}_{B, A} \| X_{t+2} - [B, A] D_{t+1} \|_F.
    \]
    \State Perform one-step policy gradient descent
    \begin{equation}\label{equ:mbgd}
    K_{t+1} = K_t - \eta \nabla \hat{C}_{t+1}(K_t),    
    \end{equation} 
    where \(\nabla \hat{C}_{t+1}(K_t)\) is the policy gradient with the estimated model \((\hat{A}_{t+1}, \hat{B}_{t+1})\).
\EndFor
\end{algorithmic}
\end{algorithm}

A major challenge is that least square identification (\ref{OLS}) becomes fragile when dynamics vary with time. More precisely, the data collected after a switch contains samples generated by different $(A_i,B_i)$, and performing ordinary least squares on such mixed–mode data returns an estimated model that may not correspond to the either the current or previous dynamics. 
  
\subsection{Closed-loop stability analysis}
For the remainder of the paper, we require a number of assumptions. First, we require persistent excitation of the state and input signals.
\begin{assum} \label{ass1}
    The sliding window data matrix has full row rank, i.e., $\textrm{rank}(D_t)=m+n$ for all $t \geq t_0$.
\end{assum}

This assumption can be ensured with probability 1 by letting $L \geq n+m$ and applying random probing noise to the control inputs
\begin{equation}
    u_t = Kx_t + e_t.
\end{equation}

We also assume a lower bound on the \textit{dwell time}, that is, how long the system remains in one mode before switching to another.
\begin{assum} \label{ass2}
    The dwell time is lower bounded by the sliding window length, i.e.,  $T_{i+1} - T_i \geq L$ for any $i \in \mathbb{N}$. 
\end{assum}

This assumption ensures that the data-series $X_t, U_t, X_{t+1}$ contains data from at most two consecutive dynamics $(A_i, B_i)$ and $(A_{i+1}, B_{i+1})$, which will ease the subsequent analysis. 

We first quantify the identification error with ordinary least-squares and sliding window data. Let the variation of the dynamics at the $i$-th switch be 
$$
\Delta_i = [B_{i+1}, A_{i+1}] - [B_{i}, A_{i}].
$$
Then, we have the following results.
\begin{lemma}
\label{error_LS}
Consider the ordinary least-squares problem with sliding window data 
\begin{equation}\label{equ:ls}
    [\hat{B}_{t}, \hat{A}_{t}] = \underset{[B, A]}{\text{arg min}} \| X_{t+1} - [B,A]D_t
    \|_F.
\end{equation} 
Then, for any $i \in \mathbb{N}$, we have
\begin{equation*}
\begin{aligned}
    \|[\hat{B}_{t}, \hat{A}_{t}] - [B_{i}, A_{i}]\| \leq \|\Delta_i\|,  &\quad T_i \leq t < T_i + L, \\
    \|[\hat{B}_{t}, \hat{A}_{t}] - [B_{i+1}, A_{i+1}]\| \leq \|\Delta_i\|,   &\quad T_i \leq t < T_i + L, \\
     [\hat{B}_{t}, \hat{A}_{t}] = [B_{i+1}, A_{i+1}], &\quad T_i + L \leq t < T_{i+1}.
\end{aligned}
\end{equation*}

\end{lemma} \vspace{2mm}

Lemma \ref{error_LS} shows that the identification error at the $i^{th}$ and $(i+1)^{th}$ mode is upper bounded by the variation of the dynamics $\|\Delta_i\|$ in $T_i \leq t < T_i + L$. For time $T_i + L \leq t < T_{i+1}$, the sliding window data is generated only by $(A_{i+1}, B_{i+1})$, and hence the identification error will be zero. By leveraging Lemma \ref{error_LS}, we can quantify the impact of using the policy gradient with the identified model instead of the true model in the policy update \eqref{equ:mbgd}.

Our main result provides conditions under which we can guarantee sequential strong stability of the closed-loop between system~\eqref{equ:switch} and Algorithm~\ref{algo2_switch}.
\begin{theorem}[\textbf{Stability in $N$ switches}]
\label{theorem_multiple_switch}
Let Assumptions \ref{ass1} and \ref{ass2} hold.
Consider $N$ switches between $N+1$ different systems. There exist constants $\nu_1$ and $\nu_{2,i}$ depending on $(A_i, B_i, Q, R)$ such that, if $\eta \leq \nu_1$, if the initial policy $K_{t_0}$ stabilizes the initial system $(A_0, B_0)$, and if
\begin{equation}
\label{conditions_multiple}
    \|\Delta_i\| \leq \nu_{2,i}, \quad i\in \{1,..., N\},
\end{equation}
then Algorithm \ref{algo2_switch} is feasible and the sequence $K_t$ is ($\overline{\kappa}, \overline{\alpha}$)-strongly stable with some constants $(\overline{\kappa}, \overline{\alpha})$. Moreover, the state is bounded by    
\begin{equation}
    \|x_t\| \leq \overline{\kappa} (1-\frac{\overline{\alpha}}{2})^{t-t_0}\|x_{t_0}\| +\frac{2\overline{\kappa}}{\overline{\alpha}} \underset{t_0\leq j < t}{\max} \|B_j e_j\|.
\end{equation}
\end{theorem} \vspace{2mm}

Theorem \ref{theorem_multiple_switch} implies that, if the stepsize  and the variation of the dynamics are sufficiently small, then the closed-loop system is sequentially strongly stable, and the state sequence is bounded above. Note that the state bound comprises two terms; the first decreases exponentially with respect to the initial state after an overshoot with magnitude $\bar{\kappa}$, and the second uniformly bounds the bias induced by the probing noise.

Recall Assumption \ref{ass2}, which imposes that for each $i$, we have $T_{i+1}-T_i \geq L$, to make sure that the sliding window data is generated by at most two systems. Therefore, condition (\ref{conditions_multiple}) implicitly yields a bound on the variation of the system matrices.  

We can now apply Theorem~\ref{theorem_multiple_switch} repeatedly, ensuring that the algorithm has enough time to converge toward mode-specific optimality before switching occurs. As a result, the system’s cost function remains bounded at all times. In turn, this allows us to provide a stricter condition on the dwell time, under which Algorithm \ref{algo2_switch} guarantees sequential stability for an unbounded number of switches.

\begin{coro}[\textbf{Infinite switches}]
\label{thm:dwell}
Let Assumption \ref{ass1} hold. Suppose that there exist variables $\nu_{1,i}$ and $\nu_{2,i}$ depending on $(A_i, B_i, Q, R)$ such that the initial policy $K_{t_0}$ stabilizes the initial system $(A_0, B_0)$, $\eta \leq \nu_{1,i}$, condition (\ref{conditions_multiple}) holds, and for all $i$
\begin{equation}\label{equ:dwell_time}
T_{i+1}-T_i \geq  \frac{2\mu_i}{\eta}\log\left( \frac{1}{\overline{C}_i - C_i^*} \right) +L
\end{equation}
then Algorithm \ref{algo2_switch} remains feasible for an infinite number of switches. Moreover, the state is bounded by    
\begin{equation}\label{eq:statebound}
    \|x_t\| \leq \nu_{3,i} (1-\frac{\nu_{4,i}}{2})^{t-t_0}\|x_{t_0}\| +\frac{2\nu_{5,i}}{\nu_{4,i}} \underset{t_0\leq j < t}{\max} \|B_j e_j\|,
\end{equation} 
where $\nu_{3,i}, \nu_{5_i} > 1$ and $0<\nu_{4,i}<1$. 
\end{coro} \vspace{2mm}
 
In the following, we discuss limitations and extensions of this result.

\textbf{Uniform bounds.} Note that conditions \eqref{conditions_multiple} and \eqref{equ:dwell_time} are mode-dependent and therefore dependent on time. Specifically, each variation $\Delta_i$ and dwell-time is determined by the corresponding mode parameters $(A_i, B_i, C_i^*, \overline{C_i})$. For Corollary~\ref{thm:dwell} to yield a bound on the states that is independent of time, we can assume a \textit{common} upper bound $\overline{C} \geq C_i(K_t)$ for all $i$. In turn, additional boundedness assumptions on the system matrices can be used to derive such uniform upper bounds on the cost. 

\textbf{Switched systems with process noise.} While process noise is not explicitly considered in our analysis, it is straightforward to adapt the results for the time invariant case \cite{zhao2025policy} to this problem. In particular, this requires changing $\nu_{2,i}$ in \eqref{conditions_multiple} to denote a \textit{robust} margin of stability and adding an additional term in \eqref{eq:statebound} to account for the noise. 

\textbf{Reducing dwell time with multiple-step policy gradient.} In Algorithm \ref{algo2_switch}, we perform only a single gradient descent step at each time instance. If online computation allows, we can instead conduct multiple steps of gradient descent per time, increasing the rate of convergence. This means that the algorithm can rapidly converge to the optimal policy for the varying dynamics and the dwell time can be reduced.

\textbf{Relaxing the dwell time.} For ease of exposition, we employ Assumption~\ref{ass2} and consider switched systems \eqref{equ:switch} with a bound on the dwell time. Nevertheless, results similar to Lemma~\ref{error_LS} can be derived  with systems that switch more often. Assuming an \textit{average} dwell time, one can then derive sequential strong stability in expectation. In this case, considering data with a forgetting factor might be more natural than using sliding window data \eqref{equ:swd}.


\section{Simulations}
We run a simulation to show the the upper bound of the LQR cost and the stability of the state. We choose the following parameters: 
\[ A_0 \!\!=\! \begin{bmatrix}
        -0.13 & 0.14 & \!\!-0.29 & 0.28 \\
        0.48 & 0.09 & 0.41 & 0.30 \\
        -0.01 & 0.04 & 0.17 & 0.43 \\
        0.14 & 0.31 & \!\!-0.29 & \!\!-0.10
        \end{bmatrix}\!\!, B_0 = \begin{bmatrix}
        1.63 & 0.93 \\
        0.26 & 1.79 \\
        1.46 & 1.18 \\
        0.77 & 0.11 
    \end{bmatrix}.\]
Moreover, we assume that the dynamics follow a random walk, that is, for each $i$, we take $\bar{A_i},\bar{B}_i$ matrices for which each element is drawn from the standard normal distribution, let
$A_{i+1} = A_i + 0.1\bar{A}_i$, and $B_{i+1} = B_i + 0.1\bar{B}_i.$
Lastly, for the LQR and algorithm parameters, we take $Q$ and $R$ as identity matrices, and
$L = 25, T_{i+1} - T_i = 30, \Delta t_i = 5.$

\subsection{Upper boundedness of the LQR cost}
\begin{figure}[t]
  \centering
\includegraphics[width=1\linewidth]{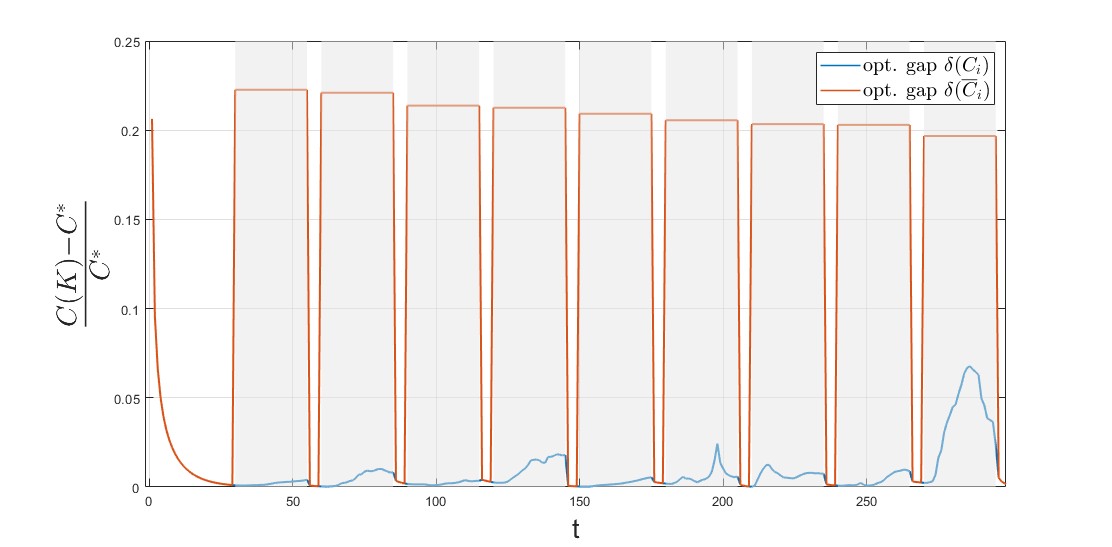}
  \caption{Comparison between real cost function $C_i(K_t)$ and its upper bound $\overline{C}_i$}
  \label{fig:opt_gap_bounds}
\end{figure} 
As established in Section III, the cost satisfies the inequality $C_{i+1}(K_t) \leq \overline{C}_{i+1}$ during the transition period following the $i^{th}$ switch, i.e. $t \in [T_i+1, T_i + L]$. In the simulations, we illustrate the evolution of the cost function and its corresponding upper bound in terms of the optimality gap, defined as
\begin{equation}
    \delta(C(K)) = \frac{C(K) - C^*}{C^*}.
\end{equation}

In Figure~\ref{fig:opt_gap_bounds}, the upper bound of the cost was computed as stated in Lemma \ref{lemma_c1_upperbound} in the Appendix. The gray areas represent the transition periods $t \in [T_i+1, T_i+L]$, where data series spans two consecutive modes. Figure~\ref{fig:opt_gap_bounds} shows that the cost upper bound is satisfied throughout each transition period. It also indicates that the cost function quickly returns to near-optimal values after each transition, suggesting that a significantly shorter dwell time could be feasible.

\subsection{State norm}
We also show the upper bound of the state norm with our setting. 
\begin{figure}[t]
  \centering
\includegraphics[width=1\linewidth]{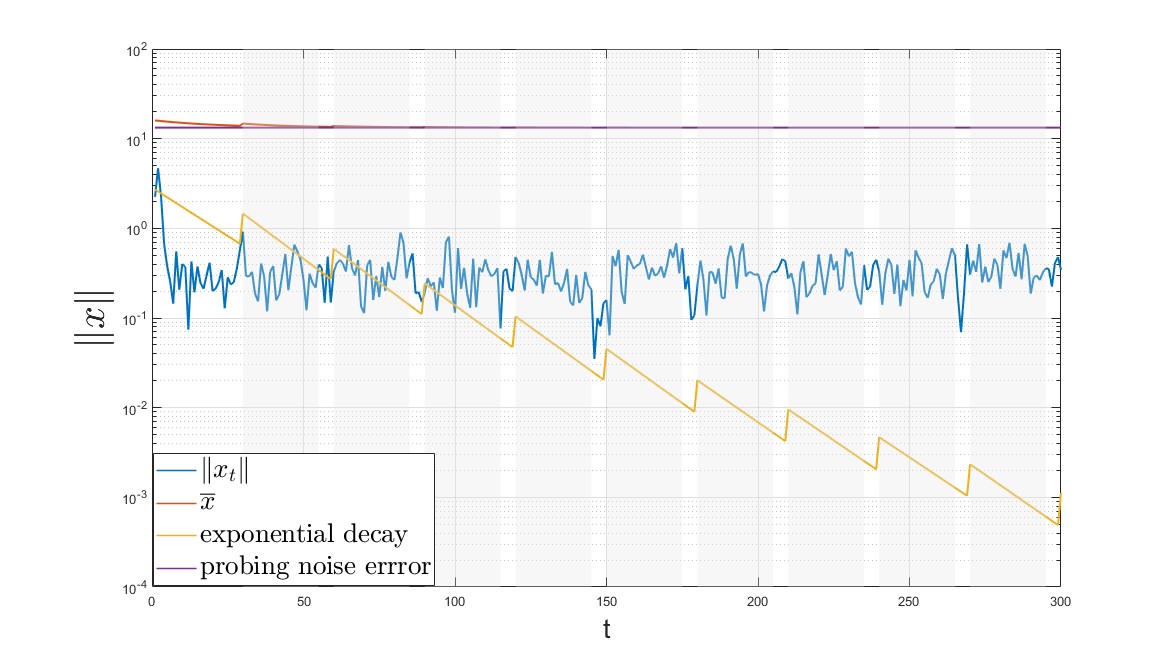}
  \caption{Comparison between the state $\|x_t\|$ norm and its upper bound.}
  \label{fig:state_norm_bounds}
\end{figure} 
As illustrated in Figure \ref{fig:state_norm_bounds}, the state norm remains consistently small and bounded across several switching instances. A change in the upper bound is observed at each switch, reflecting the time-varying nature of the variables $\nu_3$, $\nu_4$, and $\nu_5$ in Corollary~\ref{thm:dwell}. However, the contribution of the probing-noise-induced error  term (purple line) quickly dominates the exponential decay (yellow line), making our upper bound too conservative.
\section{Conclusions}
In this paper, we have proposed a policy gradient adaptive control method to stabilize unknown linear switched systems, where the system changes after a dwell time. We have shown closed-loop stability of the proposed method and demonstrated our results via simulations.

In future, it would be interesting to consider continuously varying systems without dwell time. Moreover, adopting discounted LQR cost as the performance index and analyzing the performance are worth investigating.

\appendices
\section{Proofs}
\subsection{Proof of Lemma \ref{error_LS} (Identification Error)} 
Define \(W_t = \Delta_iD_t S\), where \(S=\mathrm{diag}(1,\ldots,1,0,\ldots,0)\) selects the pre-switch columns. Then, we formulate the subspace relation as 
\begin{equation}
    X_{t+1} = A_i X_t +B_iU_t + W_t.
\end{equation}
 For the ordinary least-squares problem \eqref{equ:ls}, it holds that $\|[\hat{B}_{t+1}, \hat{A}_{t+1}] - [B_i, A_i]\| = \|\Delta_i\, D_t S D_t^\dagger\| \leq \|\Delta_i\|\,\|D_t S D_t^\dagger\|$. Using the definition of induced matrix norm and letting \(M=D_t^\top D_t\succ0\), it follows that
\[
\|D_t S D_t^\dagger\|^2
\leq \sup_{z\neq 0} \frac{z^\top S M S z}{z^\top M z}
\leq 1,
\]
where the last inequality follows from \(SMS\preceq M\). Hence, it holds that \(\|D_t S D_t^\dagger\|\le 1\), which completes the proof.

\subsection{Proof of upper boundedness of the LQR cost}

Unlike in \cite{zhao2025policy} for LTI systems, it is clear that analyzing the convergence of the cost function $C(K)$ is not relevant as we assume unknown switches in the dynamics which will induce non-optimality in the cost. A more coherent approach would be to guarantee an upper bound of the cost throughout the changes. 

\begin{lemma}
\label{lemma_cost_upperbound}
Let \( A_0 + B_0K_{t_0} \) be stable. Suppose that
\begin{equation}
   \|\Delta_i\| \leq \nu_{1,i} \quad \eta \leq \nu_{2},
\end{equation}
where $\eta$ is the step-size of the gradient descent step, and $\nu_1, \nu_2$ are variables dependent of $(A_i, B_i, Q, R)$. Then,
it holds that the cost function $C_t(K_t)$ is upper bounded.
\end{lemma}

This lemma ensures that the cost function of the switch system is upper bounded at every time-step.

As stated in Assumption \ref{ass2}, the period between two consecutive switches is long enough so that data series can only at most span two consecutive modes $(i, i+1)$. Thus, we can treat each switch independently. We focus on periods where data series span two modes, i.e $t \in ]T_i, T_i + L]$. We first demonstrate that $C_i(K)$, the cost function associated with pre-switch system $(A_i, B_i)$, remains upper bounded after the switch. We then demonstrated that $C_{i+1}(K)$, cost function associated with post-switch system $(A_{i+1}, B_{i+1})$, is upper bounded. 
\begin{lemma}[Lemma 11, \cite{zhao2025policy}]
\label{lemma_bounds}
For \(K\in\mathcal{S}_i\), it holds that
(i)\(\|\Sigma\| \le \dfrac{C(K)}{\sigma(Q)}\),
(ii) \(\|P\|\le C(K)\),
(iii) \(\|K\|_F \le \Big( \dfrac{C(K)}{\sigma(R)} \Big)^{1/2}\).
\end{lemma} 

\begin{lemma}[Lemma 12, \cite{zhao2025policy}]
\label{lemma_stability_margin}
Let \(A\in\mathbb{R}^{n\times n}\) be stable and \(\Sigma(A)=I + A\Sigma(A)A^\top\).
If \(\|A'-A\| \le (4\|\Sigma(A)\|(1+\|A\|))^{-1}\), then \(A'\) is stable and
\(\|\Sigma(A')-\Sigma(A)\|\le 4\|\Sigma(A)\|^2(1+\|A\|)\|A'-A\|\).
\end{lemma} \vspace{2mm}

\begin{lemma}
\label{lemma_stability}
Define
\[
p_1(a) = \frac{\underline{\sigma}(Q)}{4a\!\left(1 + \frac{a}{\underline{\sigma}(Q)}\right)\!\left(1 + \sqrt{\frac{a}{\underline{\sigma}(R)}}\right)}.
\]
If \(\|\Delta_i\| \le p_1(\overline{C_i})\), then \(\hat{A}_{t+1} + \hat{B}_{t+1}K_t\) is stable and both $\hat{C}_i(K_t), \nabla \hat{C}_i(K_t)< \infty$.
\end{lemma} \vspace{2mm}
\begin{proof}
Combining Lemma \ref{lemma_bounds} with \ref{error_LS}, it holds that
$\|(\hat{A}_{t+1} + \hat{B}_{t+1}K_t) - (A_i + B_iK_t)\| \leq \|\Delta_i\| (1 + \|K_t\|) \leq \|\Delta_i\|\left(1 + \sqrt{\frac{C_i(K_t)}{\underline{\sigma}(R)}}\right)$. 
Following with Lemma \ref{lemma_stability_margin} and applying bounds from Lemma \ref{lemma_bounds}, it must hold that
\begin{equation}
    \|\Delta_i\|(1 + \sqrt{\frac{C_i(K_t)}{\underline{\sigma}(R)}}) \leq \frac{\underline{\sigma}(Q)}{4C_i(K_t)(1 + \frac{C_i(K_t)}{\underline{\sigma}(Q)})}
\end{equation} \\
Assuming $C_i(K_t) \leq \overline{C}_i$. Consequently, $p_1(\overline{C}_i) \leq p_1(C_i(K_t))$.
\end{proof} 

We now demonstrate that $C_i(K_t) \leq \overline{C_i}$. We annotate $C_i(K)$ as  the true cost function of the system $[A_i, B_i]$ and $\hat{C}_t$ to the cost function of the estimated model $[\hat{A}_t, \hat{B}_t]$, we define $K'$ and $K''$ as
\begin{align}
    K' &= K - \eta \nabla \hat{C}_{t}(K) \\
    K'' &= K - \eta \nabla C_i(K)
\end{align}
By \cite[Lemmas 2 \& 3]{zhao2025policy}, if $\eta \leq \frac{1}{l_i}$, the gradient descent yields 
\begin{equation}
\label{convergence}
    C_i(K'') - C_i^* \leq (1 - \frac{\eta}{2\mu_0}) (C_i(K) - C_i^*)
\end{equation} 
\begin{lemma}[Lemma 14, \cite{zhao2025policy}]
\label{lemma_3_5}
Let \( K \in \mathcal{S} \). For any \( \tilde{K} \), if \( \| \tilde{K} - K \| \leq p_4 \), then \( \tilde{K} \in \mathcal{S} \) and 
\begin{align}
\| \tilde{\Sigma} - \Sigma \| &\leq p_5 \| \tilde{K} - K \| \\
| C(\tilde{K}) - C(K) | &\leq p_6 \| \tilde{K} - K \|
\end{align}
where $p_4 = \frac{\sigma_{\min}(Q) \mu}{4 C(K) \|B\| \left( \|A - BK\| + 1 \right)}$, 
$p_5 = 4 \left( \frac{C(K)}{\sigma_{\min}(Q)} \right)^2 \frac{\|B\| \left( \|A - BK\| + 1 \right)}{\mu}$ and $p_6 = p_5C(K)$.
\end{lemma} \vspace{2mm} 
\begin{lemma}
\label{lemma_gradient_error}
Let \( K \in \mathcal{S}_i \). If $\|\Delta_i\| \leq p_1$, then
\begin{equation}
\|\nabla \hat{C}_{i}(K) - \nabla C_i(K)\| \leq p_3 \|\Delta_i\|
\end{equation}
where $p_3 = \mathrm{poly}\left( \frac{C_i(K)}{\sigma(Q)}, \|A_i\|, \|B_i\|, \|R\|, \frac{1}{\sigma(R)} \right)$.
\end{lemma} \vspace{2mm}
\begin{proof} From \cite{zhao2025policy}, we have $\nabla C(K) = 2E\Sigma$, with $E = RK + B^TP(A+BK)$. Thus, $\nabla \hat{C}(K) - \nabla C(K) = 2\hat{E}\hat{\Sigma} - 2E\Sigma = 2(\hat{E} - E) \hat{\Sigma} - 2 E(\Sigma - \hat{\Sigma})$. First, we focus on the second term $2 E(\Sigma - \hat{\Sigma})$. By \cite[Lemma 11]{fazel2018global}, it follows that
\begin{equation}
\label{EE}
    \text{Tr}(E^TE) \leq \|R + B^TPB\|(C(K) - C^*).
\end{equation}
Combining Lemma \ref{lemma_stability_margin} \& \ref{lemma_stability} and defining $p_2 = \frac{C(K)}{\underline{\sigma}(Q) p_1}$, it holds 
\begin{equation}
\label{Sigma_diff}
    \|\hat{\Sigma} - \Sigma\| \leq \textit{p}_2 \|[\hat{B}, \hat{A}] - [B, A]\| \leq \textit{p}_2 \|\Delta \|.
\end{equation}
Putting (\ref{EE}) and (\ref{Sigma_diff}) together, it follows that 
\begin{equation}
    \|2 E(\Sigma - \hat{\Sigma})\| \leq 2 \|R + B^TPB\|(C(K) - C^*) \textit{p}_2 \|\Delta \|.
\end{equation}
To bound the term $2(\hat{E} - E) \hat{\Sigma}$, we note that $ E - \hat{E} = B^TP((A+BK) - (\hat{A}+\hat{B}K)) + (B^TP-\hat{B}^T\hat{P})(\hat{A}+\hat{B}K)$.
Bounding each term, we have $\|B^TP((A+BK) - (\hat{A}+\hat{B}K))\| \leq \|B^TP\|(1+\|K\|)\|\Delta \|$, and $ \|\hat{A}+\hat{B}K\| \leq \|A+BK\| +(1+\|K\|)p_1$. For the term $B^TP-\hat{B}^T\hat{P}$, it is developed as 
$B^TP-\hat{B}^T\hat{P}= \hat{B}^T (P - \hat{P})  +  (B^T - \hat{B}^T)P$. We note that $\|\hat{B}\| \leq \|B\| + p_1$, $\|B^T - \hat{B}^T\| \leq \|\Delta \|$, $\|P - \hat{P}\| \leq \|Q + K^TRK\| \|\Sigma - \hat{\Sigma}\| \leq C(K) p_2  \|\Delta \|$, $\|\hat{\Sigma}\| \leq 2\|\Sigma\|$. Therefore, $\|B^TP-\hat{B}^T\hat{P}\| \leq (\|B\| +p_1)\|Q + K^TRK\| p_2 \|\Delta \| + \|\Delta\| \|P\|$. \vspace{2mm} \\
Further, it follows that 
    $\|E - \hat{E}\| \leq \|B^TP\|(1+\|K\|)\|\Delta \| + \left((\|B\| +p_1)\|Q + K^TRK\| p_2  + \|P\| \right) \times \left(\|A+BK\| +(1+\|K\|)p_1 \right) \|\Delta\|$.
Noting that $\|P\|, \|\Sigma\|, \|K\|$ can be upper bounded by a polynomial in ($C(K)$, $\sigma(Q)$,$\frac{1}{\sigma(R)}$, $\|R\|$, $\|A\|$, $\|B\|$), the proof is completed.
\end{proof} 

This shows that the error in the gradient, consequently to the model mismatch, is bounded and proportional to the magnitude of the variation. To bound the error in the cost function induced by inaccuracies in the gradient, we make use of the following lemma.
\begin{lemma}
\label{lemma_3_7}
Let \( K \in \mathcal{S}_i \). There exists a polynomial \( p_7 \) in 
\(
\left( \frac{\sigma(Q)}{C_i(K)}, \frac{1}{\|A_i\|}, \frac{1}{\|B_i\|}, \frac{1}{\|R\|}, \sigma(R) \right)
\)
such that, if 
\begin{equation}
\|\Delta\| \leq \overline{p}_1 \quad \text{and} \quad \eta \leq \min \left\{ \frac{p_4}{p_3 \|\Delta\|}, p_7 \right\},
\end{equation}
then it holds that $|C_i(K'') - C_i(K')| \leq \eta p_3 p_6 \|\Delta_i\|.$

\end{lemma} \vspace{2mm}
Lemma~\ref{lemma_3_7}, adapted from \cite[Lemma 15]{zhao2025policy},  provides an upper bound on the error in the cost function resulting from inaccuracies in the gradient descent step. \vspace{2mm} \\ 
\begin{proof} Let $p_7$ be the polynomial of the step-size in \cite[Theorem 7]{fazel2018global}. Then, it follows from \cite[Theorem 7]{fazel2018global} that the update returns a stabilizing policy $K''$. It also holds that $\|K' - K''\| \leq \eta\|\nabla \hat{C}(K) - \nabla C(K)\| \leq \eta p_3 \|\Delta\|$. By \textit{Lemma 6} and if $\eta \leq \frac{p4}{p_3 \|\Delta\|}$, we obtain  $|C(K'') - C(K')| \leq \eta p_3 p_6 \|\Delta\|.$
\end{proof} 

Going back to \ref{convergence} and noting Lemma \ref{lemma_3_7}, we obtain
\begin{equation}
\label{Conv_residual_error}
    C_i(K') - C_i^* \leq (1 - \frac{\eta}{2\mu_i}) (C_i(K) - C_i^*) + \eta p_3 p_6 \|\Delta_i\|.
\end{equation}
Let $\overline{p}_j$ polynomials associated with $\overline{C}_i$ and $p^*_j$ polynomials associated with $C_i^*$, then;
\begin{lemma}
Let \( K_T \in \mathcal{S}_i \).
If 
$$
    \|\Delta_i\| \leq \min \left\{ \frac{1}{2\mu_i \overline{p}_3 \overline{p}_6}, \overline{p}_1 \right\}, \eta \leq \min \left\{ \frac{1}{\overline{l}_i}, \frac{p^*_4}{\overline{p}_3 \|\Delta\|}, \overline{p}_7\right\},
$$
then it holds that 
\begin{equation}
    C_i(K_{t}) \leq \max \left\{ C_i(K_{T_i}), C_i^* +1 \right \}.
\end{equation}   
\end{lemma} 
Following the convergence result in (\ref{Conv_residual_error}), Lemma~\ref{lemma_cost_upperbound} establishes an upper bound on $C_i$, provided that $\|\Delta_i\|$ and $\eta$ are sufficiently small. \vspace{2mm} \\
\begin{proof} From \eqref{convergence} and Lemma \ref{lemma_3_7}, the convergence of $C_i$ becomes 
\begin{equation*}
\label{conv_residual}
C_i(K_{t+1}) - C_i^* \leq (1 - \frac{\eta}{2\mu_i}) (C_i(K_t) - C_i^*) + \eta \overline{p}_3 \overline{p}_6 \|\Delta_i\|.
\end{equation*}
With the condition $\|\Delta_i\| \leq \frac{1}{2\mu_i \overline{p}_3 \overline{p}_6}$, it follows from (\ref{conv_residual}) that
\begin{equation}
\label{conv_bound}
C_i(K_{t+1}) - C_i^* \leq (1 - \frac{\eta}{2\mu_i}) (C_i(K_t) - C_i^*) + \frac{\eta}{2\mu_i}.
\end{equation}
Assuming \( K_t \in \mathcal{S}_i \). If $C_i(K_t) \geq C_i^* + 1$, (\ref{conv_bound}) becomes
\begin{equation}
\label{proof1}
\begin{aligned}
        C_i(K_{t+1}) - C_i(K_t) &\leq - \frac{\eta}{2\mu_i} (C_i(K_t) - C_i^*) + \frac{\eta}{2\mu_i} \\
        &\leq - \frac{\eta}{2\mu_i}  + \frac{\eta}{2\mu_i} = 0\\
\end{aligned}
\end{equation}
If $C_i(K_t) \leq C_i^* + 1$, (\ref{conv_bound}) becomes
\begin{equation}
\label{proof2}
\begin{aligned}
    C_i(K_{t+1}) - C_i^* &\leq (1- \frac{\eta}{2\mu_i}) (C_i(K_t) - C_i^*) + \frac{\eta}{2\mu_i} \\
    &\leq 1- \frac{\eta}{2\mu_i} + \frac{\eta}{2\mu_i} = 1\\
\end{aligned}
\end{equation}
As $C_i(K_{T_i}) \in \mathcal{S}_i$ and with (\ref{proof1}) and (\ref{proof2}), we obtain by induction that  
    $C_i(K_{t}) \leq \max\{C_i(K_{T_i}), C_i^* + 1\} =: \overline{C}_i$,
which completes the proof. 
\end{proof} \vspace{2mm}
Now, we demonstrate that $C_{i+1}(K_t)$ is upper bounded for $t \in ]T_i, T_i + L]$. 
\begin{lemma}
\label{lemma_c1_upperbound}
    If the condition of Lemma \ref{lemma_cost_upperbound} holds, then 
\begin{equation}
\label{3.42}
    \overline{C}_{i+1} \leq \overline{C}_i(1+\overline{p}_2 \|\Delta_i\|).
\end{equation}
\end{lemma} 
Lemma \ref{lemma_c1_upperbound} provides an upper bound for the cost of the new system $C_{i+1}$, proportional to $\overline{C}_i$ and their difference $\|\Delta_i\|$. \vspace{2mm} \\
\begin{proof}: By definition, 
$$
    |C_{i+1}(K) - C_i(K)| = \left|\text{Tr}\left( (Q+K^\top RK)(\Sigma_{i+1} - \Sigma_{i}) \right) \right|.
$$
By Lemma \ref{lemma_3_5}, we know that 
$
    \|\Sigma_{i+1} - \Sigma_i\| \leq p_2 \|\Delta_i\|.
$
Therefore, $\left|\text{Tr}\left( (Q+K^\top RK)(\Sigma_{i+1} - \Sigma_{i}) \right)\right| \leq \text{Tr}(Q+K^\top RK) \|\Sigma_{i+1} - \Sigma_{i}\|, \leq \text{Tr}(Q+K^\top RK) p_2 \|\Delta_i\|$.
From the definition of $\Sigma$,
$\underline{\sigma}(\Sigma) \geq 1$,
which implies that 
\begin{equation}
\begin{aligned}
    C(K) &\geq \text{Tr}(Q + K^\top RK) \underline{\sigma}(\Sigma), \\
    \frac{C(K)}{\underline{\sigma}(\Sigma)} &\geq \text{Tr}(Q + K^\top RK), \\
    C(K) &\geq \text{Tr}(Q + K^\top RK).
\end{aligned}
\end{equation}
As $C_i$ is an upper bound for any $K_t$, it follows that $\text{Tr}(Q + K_t^\top RK_t) \leq \overline{C}_i.$ Putting everything together, we obtain
\begin{equation}
\begin{aligned}
    |C_{i+1}(K_t) - C_i(K_t)| &\leq \text{Tr}(Q + K_t^\top RK_t) \|\Sigma_{i+1} - \Sigma_i\| \\ &\leq \overline{C}_i \overline{p}_2 \|\Delta_i\|.
\end{aligned}
\end{equation}
We conclude that 
$\overline{C}_{i+1} \leq \overline{C}_i (1+\overline{p}_2 \|\Delta_i\|)$,
which is equal to (\ref{3.42}). It completes the proof for Lemma \ref{lemma_c1_upperbound}.
\end{proof}

\subsection{Proof of Theorem 1 (State Stability)}
We demonstrate sequential strong stability in the intervals $t \in ]T_i, T_i +L]$. 

\begin{lemma}[Lemma 18, \cite{zhao2025policy}]
\label{str_stab}
The policy \( K \in \mathcal{S} \) is \((\kappa, \alpha)\)-strongly stable with
\[
\kappa = \sqrt{ \frac{C(K)}{ \min\{ \sigma(R), \sigma(Q) \} } }, \quad \alpha = 1 - \sqrt{1 - \frac{1}{\kappa^2}}.
\]    
\end{lemma} \vspace{2mm}
Using Lemma \ref{str_stab} with $\overline{C}_{i+1}$, we note that every $K_t, t \in ]T_i, T_i+L]$ is $(\overline{\kappa}_{i+1}, \overline{\alpha}_{i+1})$-strongly stable. 

Recalling the definitions of $p_i, i\in \{1,2,..,7\}$, we have
\begin{lemma}
\label{seq_str_stab}
If 
\begin{equation}
\label{conditions}
    \|\Delta_i\| \leq \left\{ \frac{1}{2\mu_i \overline{p}_3 \overline{p}_6}, \overline{p}_1 \right\}, \eta = \text{min} \left\{ \frac{\overline{\alpha}_1}{2\overline{\kappa}_1^2\overline{p}_5  \overline{p}_{8}}, \frac{1}{\overline{l}_i}, \frac{p^*_4}{\overline{p}_{8}}, \overline{p}_7 \right\},
\end{equation}
then, the policy sequence $K_t$ is \((\overline{\kappa}_{i+1}, \overline{\alpha}_{i+1})\)-strongly stable for $t \in ]T_i, T_i+L]$.   
\end{lemma}  \vspace{2mm}  
\begin{proof}
    The proof follows that of Lemma 17 in \cite{zhao2025policy}.
\end{proof}
\subsection{Convergence of the state}
\begin{lemma}
\label{conv_state}
If conditions (\ref{conditions}) hold, then sequential strong stability is guaranteed and 
$$
    \|x_t\| \leq \overline{\kappa}_{i+1} (1-\frac{\overline{\alpha}_{i+1}}{2})^{t-T_i-1}\|x_{T_i}\| +\frac{2\overline{\kappa}_{i+1}}{\overline{\alpha}_{i+1}} \max \|B_{i+1} e_t\|
$$
for $t \in ]T_i, T_i+L]$.
\end{lemma} \vspace{2mm} 
\begin{proof}
    The proof follows that of Lemma 18 in \cite{zhao2025policy}. 
\end{proof}
\begin{lemma}
\label{theorem}There exist constants $\nu_1, \nu_2, \nu_3$ depending on ($\overline{\kappa}_i, \overline{\alpha}_i$) for $i = [0,\dots,n]$ such that the sequence $K_t$ is $(\kappa, \alpha)$-strongly stable for constant $\kappa > 1$, $0 < \alpha < 1$ and $t\in [t_0, T_{i+1}]$, and the state is bounded
$$
    \|x_t\| \leq \nu_1 (1-\frac{\nu_2}{2})^{t-t_0-i}\|x_{t_0}\| +\frac{2\nu_3}{\nu_2} \underset{T_{i-1}\leq j < t}{\max} \|B_j e_j\|.
$$
\end{lemma} \vspace{2mm} 
\begin{proof}
    From \cite{zhao2025policy}, the sequence $K_t$ for $t \in [T_{i-1}+L, T_{i}]$ is \((\bar{\kappa}_i, \bar{\alpha}_i)\)-strongly stable (\((\bar{\kappa}_i, \bar{\alpha}_i)\) associated with $\overline{C}_i)$. Further, the sequence $K_t$ for $t \in [T_i+L, T_{i+1}]$ is \((\bar{\kappa}_{i+1}, \bar{\alpha}_{i+1})\)-strongly stable (\((\bar{\kappa}_{i+1}, \bar{\alpha}_{i+1})\) associated with $\overline{C}_{i+1})$.
\begin{equation}
\label{sss}
    \begin{aligned}
        &\|x_t\| \leq \bar{\kappa}_i (1 - \bar{\alpha}_i/2)^{t-T_{i-1}-1} \|x_{T_{i-1}}\| \\ &+ \bar{\kappa}_i \sum_{j=T_{i-1}}^{t-1} (1 - \bar{\alpha}_i/2)^{t-j-1} \|B_ie_j\|, \quad t \leq T_i \\
        &\|x_t\| \leq \bar{\kappa}_{i+1}\vspace{-2em} (1 - \bar{\alpha}_{i+1}/2)^{t-(T_i+L)-1} \|x_{T_i+L}\| \\ & + \bar{\kappa}_{i+1} \sum_{j=T_i+L}^{t-1} (1 - \bar{\alpha}_{i+1}/2)^{t-j-1} \|B_{i+1}e_j\|, t > T_i + L
    \end{aligned}
\end{equation}
Putting (\ref{sss}) and Lemma \ref{conv_state} together, we obtain
\begin{align*}
    \|x_t\| &\leq \overline{\kappa}_i (1-\frac{\overline{\alpha}_i}{2})^{t-1}\|x_{T_{i-1}}\| \\ &+\frac{2\overline{\kappa}_i}{\overline{\alpha}_i} \underset{T_{i-1}\leq j < T_i}{\max} \|B_i e_i\|, t \in ]T_{i-1}, T_i] \\
    \|x_t\| &\leq \overline{\kappa}_i \overline{\kappa}_{i+1} (1-\frac{\overline{\alpha}_i}{2})^{T_i-1}(1-\frac{\overline{\alpha}_1}{2})^{t-T_i-1}\|x_{T_{i-1}}\| \\ &+ \bar{\kappa}_i \sum_{j=T_{i-1}}^{T_i-1} (1 - \bar{\alpha}_i/2)^{T_i-j-1} \|B_ie_j\| \\ &+\bar{\kappa}_{i+1} \sum_{j=T_i}^{t-1} (1 - \bar{\alpha}_{i+1}/2)^{t-j-1} \|B_{i+1}e_j\| \\
    &\leq \overline{\kappa}_i \overline{\kappa}_{i+1} (1-\frac{\overline{\alpha}_i}{2})^{T_i-1}(1-\frac{\overline{\alpha}_{i+1}}{2})^{t-T_i-1}\|x_{T_{i-1}}\| \\ &+ \max \{\overline{\kappa}_i, \overline{\kappa}_{i+1}\} \sum_{j=T{i-1}}^{t-1} (1 - \min \{\bar{\alpha}_i, \bar{\alpha}_{i+1}\}/2)^{t-j-1} \\
    &\quad\quad\quad\quad \cdot \max\|B_{\{i,i+1\}}e_j\| \\
    \leq &\overline{\kappa}_i \overline{\kappa}_{i+1} (1-\frac{\overline{\alpha}_i}{2})^{T_i-1}(1-\frac{\overline{\alpha}_1}{2})^{t-T_i-1}\|x_{T_{i-1}}\| \\ &+ \frac{2\max \{\overline{\kappa}_i, \overline{\kappa}_{i+1}\}}{\min \{\bar{\alpha}_i, \bar{\alpha}_{i+1}\}}  \underset{T_{i-1}\leq j < t}{\max} \|B_{\{i,i+1\}}e_j\|, t \in ]T_i, T_{i+1}]
\end{align*}
Taking $\nu_1 = \overline{\kappa}_i \overline{\kappa}_{i+1}$, $\nu_2 = \min \{\overline{\alpha}_1, \overline{\alpha}_{i+1}\}$ and $\nu_3 = \max \{ \overline{\kappa}_i,  \overline{\kappa}_{i+1} \}$, it follows that 
\begin{align*}
    \|x_t\| &\leq \nu_1(1-\frac{\nu_2}{2})^{t-T_{i-1}-2}\|x_{T_{i-1}}\| \\ 
    &+\frac{2\nu_3}{\nu_2} \underset{T_{i-1}\leq j < t}{\max} \|B_{\{i,i+1\}} e_j\|.
\end{align*}

We can extend $\nu_1, \nu_2$ to $\nu_1 = \prod_{0}^{i} \overline{\kappa}_j$, $\nu_2 = \min \{\overline{\alpha}_0, \dots, \overline{\alpha}_{i}\}$ and $\nu_3 = \max \{ \kappa_0, \dots, \kappa_i\}$ for $i = [0, \dots, n]$. Thus, the state norm becomes, 
\begin{equation}
    \|x_t\| \leq \nu_1(1-\frac{\nu_2}{2})^{t-t_0-i}\|x_{t_0}\| \\ +\frac{2\nu_3}{\nu_2} \underset{0\leq j < t}{\max} \|B_j e_j\|.
\end{equation}
This proves Lemma \ref{theorem}.\end{proof}
Replacing $\overline{\kappa} = \prod_{0}^{n} \overline{\kappa}_j(1-\frac{\overline{\alpha}_j}{2})^{-1}$ and $\overline{\alpha} = \min \{\overline{\alpha}_0, \dots, \overline{\alpha}_{n}\}$, the proof of Theorem \ref{theorem_multiple_switch} is completed. 
\subsection{Proof of Corollary 1 (Dwell time condition)}
We proved that for two consecutive mode, their associated cost functions have the relation \begin{equation}
    \overline{C}_{i+1} \leq \overline{C}_i (1+\overline{p}_2 \|\Delta_i\|).
\end{equation} 
We also know the convergence rate of a noiseless LTI system (which occurs for $\forall t \in ]T_{i-1} + L, T_i]$) is 
\begin{equation}
C_i(K_{t+1}) - C_i^* \leq (1 - \frac{\eta}{2\mu_i}) (C_i(K_t) -C_i^*),
\end{equation}
and that due to identification error, our cost function is bounded by
\begin{equation}
    C_i \leq \min \{ C_i(K_{T_i}), C_i^* + 1\}.
\end{equation}
 during the period $ t \in ]T_i, T_i + L]$. To ensure a bound for the infinite horizon scenario, we impose that the cost function $C_i$ must decrease to at least $C_i^* + 1$ before a new switch. Applying this dwell time condition to our convergence equation, we have:
\begin{equation}
\begin{aligned}
    C_i(K_{T_{i}}) - C_i^* &\leq (1 - \frac{\eta}{2\mu_i})^{\Delta t_i} (C_i(K_{T_{i-1}+L}) -C_i^*) \\
     C_i(K_{T_{i}}) - C_i^* &\leq (1 - \frac{\eta}{2\mu_i})^{\Delta t_i} (\overline{C}_i -C_i^*) 
    \leq 1,
\end{aligned}
\end{equation}
which leads to 
\begin{equation}
\begin{aligned}
    \Delta t_i \geq \frac{\log\left( \frac{1}{\overline{C}_i - C_i^*} \right)}{-\log\left(1 - \frac{\eta}{2\mu_i} \right)} \approx \frac{2\mu_i}{\eta}\log\left( \frac{1}{\overline{C}_i - C_i^*} \right) 
\end{aligned}
\end{equation}
for $\frac{\eta}{2\mu_i} \ll 1$.

\bibliographystyle{IEEEtran}
\bibliography{bib/mybibfile}

\end{document}